# Fast Approximation of Optimal Perturbed Long-Duration Impulsive Transfers via Deep Neural Networks

Yue-he Zhu[*], Ya-zhong Luo[†]

National University of Defense Technology, Changsha 410073, People's Republic of China

**Abstract** The design of multitarget rendezvous missions requires a method to quickly and accurately approximate the optimal transfer between any two rendezvous targets. In this paper, a deep neural network (DNN)-based method is proposed for quickly approximating optimal perturbed long-duration impulsive transfers. This kind of transfer is divided into three types according to the variation trend of the right ascension of the ascending node (RAAN) difference between the departure body and the rendezvous target. An efficient database generation method combined with a reliable optimization approach is developed. Three regression DNNs are trained individually and applied to approximate the corresponding types of transfers. The simulation results show that the well-trained DNNs are capable of quickly estimating the optimal velocity increments with a relative error of less than 3% for all the three types of transfers. The tests on the debris chains with the total velocity increments of several thousand m/s show that the estimated results can be very close to the optimized ones with a final estimation error of less than 10 m/s.

## I.    Introduction

The design of a multitarget rendezvous mission, in which the core task is to solve a sequence-undetermined global trajectory optimization problem (GTOP) that resembles the moving-target traveling salesman problem (TSP), is usually challenging.

---

[*] Ph.D. Candidate, College of Aerospace Science and Engineering; zhuyuehe@nudt.edu.cn.

[†] Professor, College of Aerospace Science and Engineering; luoyz@nudt.edu.cn. Senior Member AIAA.

Compared with the classical moving-target TSP and its variants [1], a perturbed multiple-impulse-based GTOP such as the problem of the ninth edition of the Global Trajectory Optimization Competition (GTOC-9) [2] is more difficult to solve because achieving the optimal velocity increments of a perturbed long-duration impulsive transfer is much more complex and time consuming than calculating the distance between two cities. Huge numbers of possible transfers are required to be evaluated, and it is apparently impracticable to optimize the velocity increments for each transfer while optimizing the rendezvous sequence. Fast estimation of optimal transfer velocity increments is necessary for the design of multitarget rendezvous missions.

Due to the difficulty in achieving the optimal transfer velocity increments, few studies focused on perturbed long-duration impulsive transfers previously, and there was a lack of an efficient method that can quickly and accurately estimate the optimal velocity increments of this kind of transfer. Drawn by the problem of GTOC-9, several analytical estimation methods emerged [3-7]. These methods were developed by the participants of GTOC-9, and most of these methods estimate the optimal velocity increments by selectively adding or taking the root-sum-square of individual velocity increments required for matching the semi-major axis, eccentricity, inclination, RAAN and phase of the rendezvous target. Using only analytical calculation, the result can be quickly obtained. However, the approximating performance is usually not satisfactory. The relative estimation error can be up to 30% for some methods [6-7]. In preliminary design, the optimal transfer chain can hardly be found if applying an estimator with such a large error even if the global search ability is powerful enough. A more reliable approximation method is required for the design of multitarget rendezvous missions.

From the mathematical point of view, the estimation of the optimal transfer velocity increments is essentially a regression problem. Numerous studies have shown that machine learning-based (ML-based) methods own the superiority on solving regression problems as long as relevant learning features are well identified and training samples are sufficient [8]. The representative ML-based works that are related to this study are noteworthy. Shang and Liu [9] applied a Gaussian process regression model to evaluate the accessibility of main-belt asteroids instead of optimizing the transfer trajectories for the candidate asteroids one by one. The simulation time was greatly saved, and the estimated results were close to the numerical optimal solutions. Zhu et al. [10] expanded Shang and Liu's work to the evaluation of round-trip transfers for manned main-belt asteroids exploration mission, and the estimation error was reduced to less than 1.2%. More relatedly, Hennes et al [11] studied the fast approximation of optimal low-thrust hop between main-belt asteroids for the design of multi-asteroid visiting missions. Several state-of-the-art ML models were employed to estimate the optimal final mass of the spacecraft, and the estimation accuracy was significantly improved compared with the traditional methods. Mereta et al [12] promoted Daniel's work to a more complicated case with multiple-revolution transfers in the near earth regime and achieved similar results. The success of the above applications inspires our attempt to apply an ML-based method to estimate the optimal velocity increments of perturbed long-duration impulsive transfers.

A good learning model is required for the estimation of optimal transfer velocity increments. As an important member of the ML family, the deep neural networks (DNNs), which refer to the artificial neural networks with more than one hidden layer [13], have increasingly been popular in recent years. A DNN with an appropriate network structure and activation function is expected to have a stronger approximation ability than traditional ML models [14], especially for the complex classification and regression problems. The significant achievements of AlphaGo [15] and OpenAI [16] have attracted large amounts of attention on DNNs. The DNN-based studies that have emerged recently in aerospace fields are also noteworthy. Sánchez-Sánchez [17-18] applied the DNNs to learn the solution of inverted pendulum stabilization and optimal landing problems, so that real-time on-board trajectory planning could be achieved. Maggiori et al. [19] proposed a fully convolutional network architecture by analyzing a state-of-the-art model and solving its concerns by construction. Their overall framework showed that convolutional neural networks could be used end-to-end to process massive satellite images and provide accurate pixelwise classifications. Furfaro et al. [20] validated that a deep recurrent neural network architecture was capable of predicting the fuel-optimal thrust from sequence of states during a powered planetary descent. Peng and Bai [21] showed that a well-trained DNN could also be combined with physics-based models to improve the orbit prediction accuracy by learning space environment information from large amounts of observed data. Owing to the powerful approximation ability,

DNNs are applied as the learning model for the approximation of optimal perturbed long-duration impulsive transfers.

For better estimating the optimal velocity increments, the characteristic of this kind of transfer is first studied. It is suggested that this kind of transfer should be divided into three types, and three regression DNNs are required to estimate the optimal velocity increments individually. The most appropriate learning features and network scale (i.e., the number of nodes and hidden layers) of these regression DNNs are investigated. The comparisons between the DNN-based method and two typical analytical methods are presented, and two debris transfer chains are further applied to validate the superiority of the DNN-based method for approximating optimal perturbed long-duration impulsive transfers.

The contributions of this paper are summarized as follows.

1) An efficient optimization approach for achieving near-optimal perturbed long-duration impulsive transfers is proposed.

2) The relationship between the optimal velocity increments and the initial RAAN difference between the departure body and the rendezvous target, as well as the transfer time is analyzed. It is first found that the transfers should be divided into three types, and the optimal velocity increments should be estimated individually.

3) A DNN-based method for quickly and accurately approximating optimal perturbed long-duration impulsive transfers is developed. This method performs much better than the analytical estimation methods used in GTOC-9 and can approximate transfer chains with a final error of less than 10 m/s.

The remainder of this paper is organized as follows. Section II briefly describes the trajectory optimization model. Section III introduces the optimization approach for achieving near-optimal transfers and validates the efficiency of the approach. Section IV studies the relationship between the optimal velocity increments and the initial RAAN difference between the departure body and the rendezvous target, as well as the transfer time. Section V presents the complete process of the DNN-based method for approximating optimal transfers. Detailed simulations for determining the most appropriate learning features and network scales and the demonstration of the DNN-based method for estimating the optimal transfer velocity increments are given in Section VI. Conclusions are drawn in Section VII.

## II. Perturbed Multiple-Impulse Trajectory Optimization Model

### A. Orbital dynamics

The motion of a spacecraft flying around the Earth can be modeled as

$$\dot{\boldsymbol{r}} = \boldsymbol{v}$$
$$\dot{\boldsymbol{v}} = -\frac{\mu}{r^3}\boldsymbol{r} + \boldsymbol{a}_p \quad , \tag{1}$$

where $\boldsymbol{r}$ and $\boldsymbol{v}$ are the position and velocity in geocentric ecliptic reference frame; $\mu$ refers to the gravitational parameter; $\boldsymbol{a}_p$ denotes the perturbation acceleration caused by the factors such as nonspherical gravity and atmospheric drag. Only the secular effect of J2 perturbation is taken into account in this study. For a spacecraft that is described by means of the osculating orbital elements, the state of the spacecraft at any epoch $t$ can be computed as

$$\begin{cases} a = a_0, \quad e = e_0, \quad i = i_0 \\ \Omega = \Omega_0 - \dfrac{3}{2} J_2 (\dfrac{R_E}{p})^2 n \cos i (t - t_0) \\ \omega = \omega_0 + \dfrac{3}{4} J_2 (\dfrac{R_E}{p})^2 n (5 \cos^2 i - 1)(t - t_0) \\ M = M_0 + [n + \dfrac{3}{4} J_2 (\dfrac{R_E}{p})^2 n \sqrt{1 - e^2} (3 \cos^2 i - 1)](t - t_0) \end{cases}, \qquad (2)$$

where $J_2$ denotes the magnitude of the J2 perturbation and $R_E$ denotes the equatorial radius of the Earth; $p = a(1 - e^2)$ and $n = \sqrt{\mu / a^3}$ are the semilatus rectum and mean motion, respectively. Once the initial state $\{a_0, e_0, i_0, \Omega_0, \omega_0, M_0, t_0\}$ and the flight time $\Delta t$ are given, the osculating orbital elements at $t = t_0 + \Delta t$ is determined and the corresponding position and velocity can be directly obtained.

## B. Optimization model

There are two kinds of design variables in the multiple-impulse trajectory optimization problem [22]. The first kind refers to the maneuver times of all the impulse, which are expressed as

$$\boldsymbol{X}_1 = T_i, \quad i = 1, 2, ..., n , \qquad (3)$$

where $n$ is the total number of maneuvers. The second kind refers to the first $n$-2 impulses, which are expressed as

$$\boldsymbol{X}_2 = \Delta \boldsymbol{V}_i, \quad i = 1, 2, ..., n - 2 , \qquad (4)$$

where $\Delta \boldsymbol{V}_i = [\Delta V_{ix}, \Delta V_{iy}, \Delta V_{iz}]$ is the $i^{\text{th}}$ impulse vector. Once $\boldsymbol{X}_1$ and $\boldsymbol{X}_2$ are determined, the spacecraft can be propagated from $T_1$ to $T_{n-1}$ base on the dynamic model presented above and the last two impulses can be computed by solving a two-point boundary-value problem (TPBVP). The goal of the problem is to minimize the total velocity increments:

$$J = \min \sum_{i=1}^{n} \| \Delta \boldsymbol{V}_i \| \tag{5}$$

For a multiple-impulse rendezvous trajectory optimization problem with maximum transfer time $\Delta T_{\max}$, all the variables in Eq. (3) must be limited to the range of [0, $\Delta T_{\max}$]. To avoid dealing with the constraints of the maneuver time sequence during optimization, the variables in Eq. (3) is modified as

$$\boldsymbol{X}_1 = \eta_i, \quad i = 1, 2, ..., n , \tag{6}$$

$$\begin{cases} \eta_i = T_i / T_{i+1}, & i = 1, 2, ..., n-1 \\ \eta_n = T_n / \Delta T_{\max} \end{cases} , \tag{7}$$

where $\eta_i (i = 1, 2, ..., n)$ are all limited to the range of [0, 1]. By this modification, the constraints of the maneuver times are always satisfied and a multiple-impulse transfer can be obtained as long as the terminal position and velocity constraints are satisfied:

$$\begin{cases} \| \boldsymbol{r}_{cf} - \boldsymbol{r}_{tf} \| \leq \varepsilon_r \\ \| \boldsymbol{v}_{cf} - \boldsymbol{v}_{tf} \| \leq \varepsilon_v \end{cases} . \tag{8}$$

$\boldsymbol{r}_{cf}, \boldsymbol{v}_{cf}$ and $\boldsymbol{r}_{tf}, \boldsymbol{v}_{tf}$ are the final states of the spacecraft and the rendezvous target, respectively; $\varepsilon_r$ and $\varepsilon_v$ are the acceptable position and velocity errors.

## III.   Approach for Optimizing Perturbed Long-Duration Impulsive Transfers

This study aims to develop a DNN-based method for quickly approximating optimal perturbed long-duration impulsive transfers. Large numbers of training samples are required, and the quality of the samples, which refers to the optimality of the transfers, must be guaranteed because it significantly influences the approximating performance. In this section, the optimization approach for achieving near-optimal solutions and the validation for the efficiency of the approach are presented.

## A.  Optimization Approach

The crucial step for achieving a constraint-satisfied multiple-impulse transfer is to solve the TPBVP between the last two impulses. Solving a two-body TPBVP is very easy by using either a classical Lambert algorithm [23] or a multiple-revolution Lambert algorithm [24]. However, it is difficult to solve the perturbed TPBVP directly. In this study, a two-step approach for optimizing perturbed long-duration impulsive transfers is applied. The process of the approach is illustrated in Figure 1. Firstly, the orbital dynamics between the last two impulses are replaced with the two-body model and a classical Lambert algorithm is used to solve the two-body TPBVP. Note that the two-body model is also used for the propagation of rendezvous target from $T_{n-1}$ to $T_n$. An improved differential evolution (DE) algorithm [25] is applied to obtain the initial solution that consists of $n$-1 perturbed transfer legs and a two-body transfer leg. Then the orbital dynamics between the last two impulses are returned to the perturbed model and a sequential quadratic programming (SQP) algorithm follows to achieve the perturbed solution based on the initial solution. The transfer leg between the last two impulses is limited to one revolution to avoid the unsuccessful convergence of the SQP caused by the excessive terminal difference between the initial solution and the perturbed one.

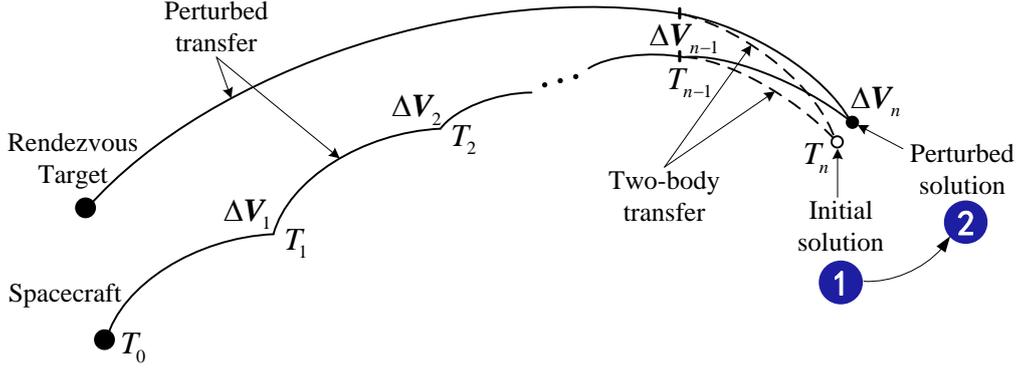

**Fig. 1 Two-step approach for obtaining a perturbed multiple-impulse rendezvous solution**

Prussing [26] pointed out that the maximum impulse number required for the optimal coplanar rendezvous problem is four. The impulse number is set to five for the non-coplanar problems considered in this study since numerous simulations show that further enlarging the impulse number contributes little to the decrease in the total velocity increments while greatly increases the optimization difficulty.

Due to the interference of huge numbers of local optima and the stochasticity of the evolutionary algorithm, the optimality of a transfer cannot be guaranteed by running the optimization algorithm only one time. 100 independent runs are implemented for each case and the best solution of the 100 runs is determined as the optimal solution.

## B. Validation for the efficiency of the approach

In the two-step approach, the transfer leg between the last two impulses is limited to one revolution. Note that the original search space is actually narrowed, and the optimal solution may be missed with such a limitation. The influence of the limitation on achieving optimal long-duration impulsive transfers is investigated and the quality of the best solution that is selected from 100 independent runs is verified.

Yang et al. [27] developed a method using homotopic perturbed Lambert algorithm to solve this problem. By introducing a set of middle target points along the position offset vector and a new targeting technique, the homotopy-based Lambert algorithm is capable of solving perturbed long-duration TPBVP, and the perturbed solution can be achieved directly without the help of SQP. Due to the lack of the method for validating the optimality, the homotopy-based approach is applied to make a comparison with the two-step approach. For fairness, the improved DE algorithm [25] is also applied as the optimizer for the homotopy-based approach.

**Table 1 Ten transfer cases for testing the efficiency of the two-step approach**

| Case | Classic orbital elements, (m, -, deg, deg, deg, deg) | | Maximum transfer time, day |
| --- | --- | --- | --- |
| | Departure body | Rendezvous target | |
| 1 | [7102019.008, 0.0033, 98.173, 1.258, 188.448, 221.186] | [7113158.741, 0.0133, 98.524, 0, 164.332, 65.562] | 9.724 |
| 2 | [7042245.022, 0.0059, 99.277, 5.268, 252.645, 259.213] | [7004095.428, 0.0151, 100.634, 0, 298.769, 222.769] | 19.989 |
| 3 | [6996726.169, 0.0138, 99.902, 0, 359.571, 135.360] | [7064707.883, 0.0199, 96.367, 15.449, 305.414, 246.734] | 23.463 |
| 4 | [6995199.808, 0.0081, 99.515, 0, 270.585, 350.836] | [7239184.088, 0.0096, 97.791, 7.218, 32.369, 200.043] | 22.202 |
| 5 | [7280713.695, 0.0028, 100.420, 0, 264.779, 228.623] | [7110280.222, 0.0134, 96.592, 9.866, 286.569, 273.241] | 28.049 |
| 6 | [7207996.616, 0.0007, 99.455, 0, 320.522, 22.026] | [7283641.352, 0.0049, 98.135, 2.725, 337.486, 227.488] | 15.081 |
| 7 | [7147223.604, 0.0067, 99.207, 0, 276.253, 27.814] | [7008415.946, 0.0063, 96.439, 4.885, 334.205, 20.624] | 18.306 |
| 8 | [6963079.862, 0.0194, 97.666, 0, 184.209, 319.741] | [6933400.666, 0.0028, 97.900, 3.691, 288.044, 283.661] | 2.582 |
| 9 | [7151921.596, 0.0156, 96.485, 0, 20.478, 43.734] | [6977324.893, 0.0067, 96.184, 3.776, 233.778, 276.058] | 5.481 |
| 10 | [6991377.623, 0.0027, 99.166, 4.705, 50.878, 200.219] | [7132670.622, 0.0199, 96.617, 0, 88.349, 171.312] | 16.876 |

Ten transfer cases are applied for the comparison, and the details are given in Table 1. 100 independent runs are implemented for both of the two approaches and the results are listed in Table 2. We can see from Table 2 that the minima of the two-step approach are almost the same as the ones of the homotopy-based approach for all the tens cases. The results indicate that limiting the last transfer leg to the last one revolution influences little on solving a long-duration five-impulse rendezvous problem. Even if the true optimal solution may be missed with such a limitation, the best solution in 100 independent runs can still be very close to the true optimal one.

**Table 2 Optimization results of the two approaches on the test cases**

| Case | Approach | Maximum, m/s | Minimum, m/s | Average time, s | Case | Approach | Maximum, m/s | Minimum, m/s | Average time, s |
|------|----------|--------------|--------------|-----------------|------|----------|--------------|--------------|-----------------|
| 1 | Two-step | 229.02 | 119.46 | 5.48 | 6 | Two-step | 195.71 | 177.70 | 5.32 |
| | Homotopy | 138.06 | 119.13 | 243.42 | | Homotopy | 185.19 | 177.70 | 337.74 |
| 2 | Two-step | 299.04 | 244.61 | 5.89 | 7 | Two-step | 407.01 | 369.89 | 5.49 |
| | Homotopy | 265.68 | 243.90 | 369.92 | | Homotopy | 379.54 | 369.91 | 388.52 |
| 3 | Two-step | 823.14 | 627.36 | 5.91 | 8 | Two-step | 695.37 | 496.21 | 4.88 |
| | Homotopy | 667.39 | 626.43 | 402.43 | | Homotopy | 534.68 | 494.60 | 144.41 |
| 4 | Two-step | 321.97 | 256.99 | 5.90 | 9 | Two-step | 831.02 | 646.80 | 5.11 |
| | Homotopy | 270.56 | 256.92 | 473.65 | | Homotopy | 729.23 | 645.70 | 194.62 |
| 5 | Two-step | 574.99 | 504.03 | 6.28 | 10 | Two-step | 588.56 | 479.96 | 5.29 |
| | Homotopy | 515.80 | 502.88 | 623.29 | | Homotopy | 528.92 | 478.44 | 354.30 |

Moreover, Table 2 shows that the maxima of the homotopy-based approach are all smaller than that of the two-step approach. This results indicate that the homotopy-based approach is more stable, and a solution obtained by the homotopy-based approach is expected to be closer to the optimum if running the optimization algorithm only one time. However, it is evident that much more time is required for the homotopy-based approach. In this study, the efficiency of achieving a near-optimal solution through many runs but not the optimization stability is more important because neither of the two approaches can guarantee to the optimum in one run and multiple runs are necessary to generate a near-optimal solution. From this point of view, the two-step approach is more efficient and thus be applied.

# IV.  Characteristic Analyses for Optimal Transfers

A multiple-impulse transfer can be essentially seen as a process to eliminate the semi-major axis difference, eccentricity difference, orbital plane angle difference and phase difference between the departure body and the rendezvous target through several maneuvers. The optimal velocity increments of a long-duration impulsive transfer are mainly determined by the first three differences because the large number of revolutions allows for the adjustment of the phase difference. For two-body cases, the optimal transfer velocity increments are almost constant because the first three differences are all fixed. However, it is quite different for the perturbed ones. The RAANs of the departure body and rendezvous target are no longer constant, and the orbital plane angle difference varies with the increase in the transfer time. This can result in the possible variation of the optimal transfer velocity increments. For better approximating optimal transfers, the relationship between the optimal velocity increments and the initial RAAN difference between the departure body and the rendezvous target, as well as the transfer time is investigated in this section.

The transfers between the bodies with different RAAN variation rates are considered. These transfers can firstly be divided into two kinds according to the variation trend of the RAAN difference between the departure body and the rendezvous target. The first kind refers to the transfers of which the RAAN difference tends to become smaller and smaller, and the second kind refers to the ones of that with opposite trend. Eight groups of transfer cases are tested, where the departure body and the rendezvous target are the same for all the cases in the eight groups. The transfers in the same group share the same initial states of the departure body and rendezvous target but with different maximum transfer times. The initial orbit elements of the departure body and rendezvous target for the transfer cases in the first group are listed in Table 3. Starting from the RAAN difference of -4 deg, the departure body and the rendezvous target are propagated forward based on the perturbed dynamic model, and the initial RAAN differences between the departure body and the rendezvous target are set to -3 deg, -2 deg, -1 deg, 1 deg, 2 deg, 3 deg, 4 deg, respectively, for the transfer cases in the other seven groups. The maximum transfer time increases from one day to tens of or even hundreds of days. All the transfer cases are obtained based on the optimization approach presented in Sec. III.

**Table 3 Initial orbit elements of the departure body and rendezvous target for the transfer cases**

|  | $a$, m | $e$, - | $i$, deg | $\Omega$, deg | $\omega$, deg | $f$, deg |
|---|---|---|---|---|---|---|
| Departure body | 7142116.504 | 0.006172 | 98.581 | 96 | 257.367 | 135.368 |
| Rendezvous target | 7052562.111 | 0.007721 | 97.203 | 100 | 13.265 | 311.656 |

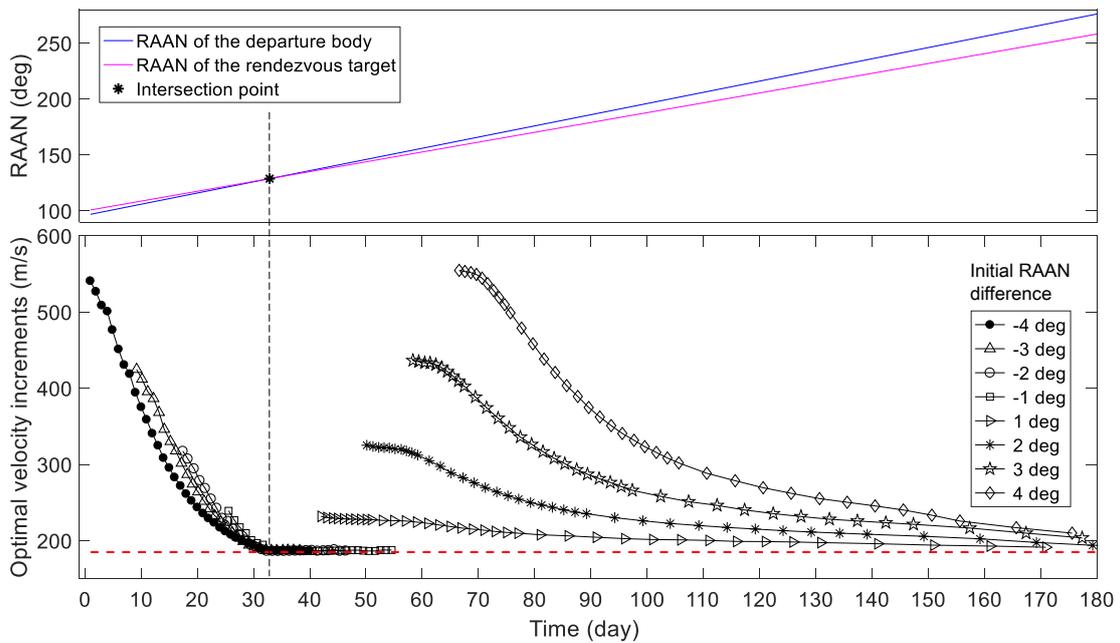

**Fig. 2 Optimal velocity increments of the transfer cases with different transfer times and initial RAAN differences between the departure body and the rendezvous target**

Figure 2 illustrates the optimal velocity increments of all the transfer cases in the eight groups. Each point corresponds to a transfer case (e.g., the third triangle corresponds to the case with an initial RAAN difference of -3 deg and a maximum transfer time of 3 days, and the fourth square corresponds to that with an initial RAAN difference of -1 deg and a maximum transfer time of 4 days).

Three interesting conclusions can be drawn from Figure 2.

1) There is a lower bound of the optimal transfer velocity increments for a transfer between any two given bodies. The lower bound is determined by the semi-major axes, eccentricities and inclinations of the departure body and rendezvous target.

2) For a transfer that is belonging to the first kind, supposing $dt$ is the free flight time required for the RAAN of the departure body to catch up with that of the rendezvous target, the optimal transfer velocity increments always decrease with the increase in the maximum transfer time if the maximum transfer time is less than $dt$ and can reach the lower bound when the maximum transfer time is approximately equal to $dt$. Further enlarging the maximum transfer time contributes little to the decrease in the optimal velocity increments if the maximum transfer time is greater than $dt$.

3) For a transfer that is belonging to the second kind, the optimal velocity increments always decrease with the increase in the maximum transfer time and can finally approach to the lower bound as long as the transfer time is sufficient.

Some other transfer cases with different initial orbit elements are also tested, and the results are similar to the situation in Figure 2.

## V.    DNN-Based Method for Approximating Optimal Transfers

Although it is impossible to analytically calculate the optimal velocity increments of a perturbed long-duration impulsive transfer, the similar variation in Sec. IV suggests that the optimal transfer velocity increments are expected to be quickly estimated with a small error using a learning-based method. The DNN-based method for estimating the optimal transfer velocity increments is presented in this section.

### A.    Implementation process

The second conclusion in Sec. IV suggests that the transfers that are belonging to the first kind can be further divided into two types according to whether the RAANs of the departure body and rendezvous target cross over the intersection point. This is because the transfers that are belonging to each of these two types show an obviously different variation trend of the optimal velocity increments. For better estimating the optimal transfer velocity increments, the transfers are ultimately divided into three types, which are named as "RAAN-closing" transfers, "RAAN-intersecting" transfers and "RAAN-separating" transfers, respectively, and the estimation of the optimal velocity increments is implemented individually. The transfers that are belonging to each of the three types are illustrated in Figure 3.

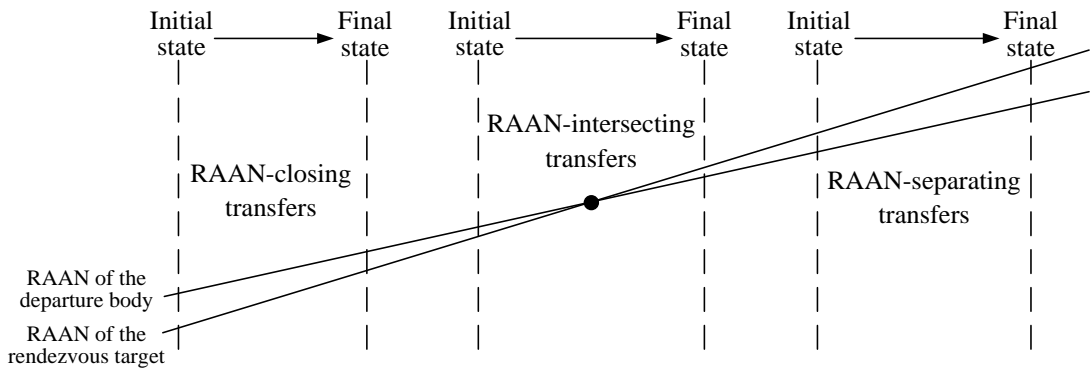

**Fig. 3 Three types of perturbed long-duration impulsive transfers**

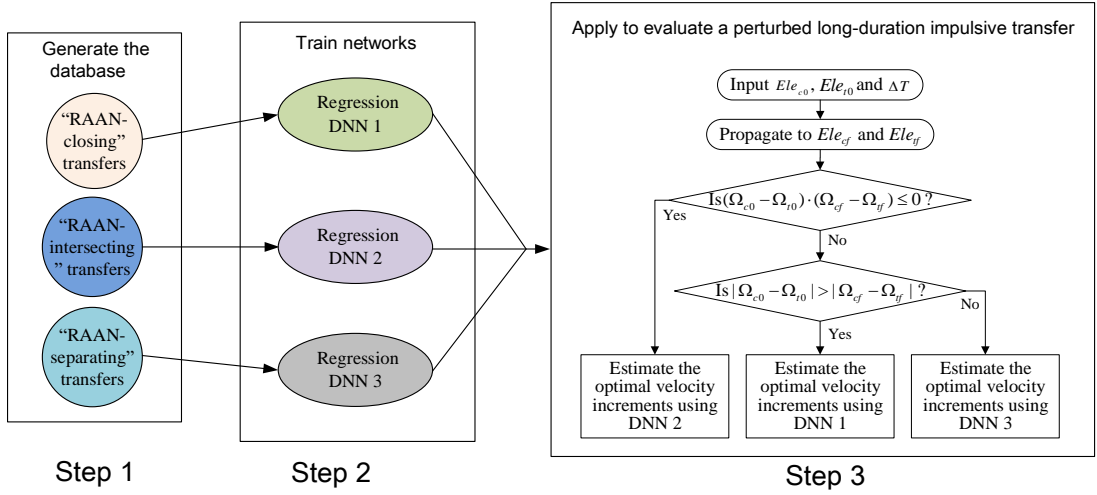

**Fig. 4 Implementation process of the DNN-based method for approximating optimal perturbed long-duration impulsive transfers**

The complete process of the DNN-based method for approximating optimal perturbed long-duration impulsive transfers is divided into three steps, which are illustrated in Figure 4.

The first step is to generate the database that contains three types of transfers. The samples should be generated according to the working conditions and parameter configurations for different problems.

The second step is to train three regression DNNs based on the corresponding samples. Note that the learning features and network scales of the three DNNs are determined in advance. The most appropriate learning features and network scales for this problem are investigated in this study, and the results can be directly applied to similar cases when training DNNs based on the new database.

The third step is to apply the well-trained DNNs to the mission design. The flow chart for approximating the perturbed long-duration impulsive transfer between any two bodies is presented in the third step, where $Ele_{c0}$, $Ele_{t0}$, $Ele_{cf}$, $Ele_{tf}$ are the initial and final orbit elements of the departure body and rendezvous target, respectively, and $\Omega_{c0}$, $\Omega_{t0}$, $\Omega_{cf}$, $\Omega_{tf}$ are the initial and final RAANs of the departure body and rendezvous target, respectively.

## B. Database generation method

Few "RAAN-intersecting" transfers can be obtained if randomly selecting two real-world central bodies and determining the transfer time when generating training samples. To improve the efficiency of obtaining "RAAN-intersecting" transfers and balance the proportions of the three types of transfers, a more efficient method is applied.

Algorithms 1 to 3 present the method of generating training samples for the three types, where $\dot{\Omega}_c$ and $\dot{\Omega}_t$ are the RAAN variation rates of the departure body and rendezvous target, and $d_1$ and $d_2$ are two parameters to determine the maxima of final and initial RAAN differences. This method first generates a training sample by randomly producing the orbit elements of the departure body and rendezvous target and then adjust the RAAN of the rendezvous target according to that of the departure body and the transfer time. Note that the samples are not the transfers between two real-world central bodies but the ones between two virtual bodies. In fact, there is no need to use the real-world central bodies for database generation because the DNN models trained by the transfers between virtual bodies can also be applied to approximate the optimal transfers between real-world central bodies as long as they have the same parameter configuration.

| Algorithm 1 Pseudocode of the process to generate an "RAAN-closing" transfer |
| --- |
| **1:** Randomly produce $Ele_{cf}$, $Ele_{tf}$ and $\Delta T$ within each range determined by the working situation of the problem |
| **2:** Calculate $\dot{\Omega}_c$ and $\dot{\Omega}_t$ |
| **3:** Randomly produce a final RAAN difference $\Delta\Omega_f$ ($\Delta\Omega_f \in (0, \boldsymbol{d}_1]$) |

| | |
|---|---|
| **4:** | **if** $\dot{\Omega}_c < \dot{\Omega}_t$ |
| | Reset $\Omega_{tf} = \Omega_{cf} - \Delta\Omega_f$ |
| | **else** |
| | Reset $\Omega_{tf} = \Omega_{cf} + \Delta\Omega_f$ |
| | **end if** |
| **5:** | Inversely propagate $Ele_{cf}, Ele_{tf}$ to $Ele_{c0}, Ele_{t0}$ (the propagation time is $\Delta T$ ) |
| **6:** | Optimize this sample ( $Ele_{c0}, Ele_{t0}$ and $\Delta T$ ) and put into the database pool 1 |
| Algorithm 2 Pseudocode of the process to generate an "RAAN-intersecting" transfer | |
| **1:** | Randomly produce $Ele_{cm}$, $Ele_{tm}$ and $\Delta T$ within each range determined by the working situation of the problem |
| **2:** | Calculate $\dot{\Omega}_c$ and $\dot{\Omega}_t$ |
| **3:** | Randomly produce the intersecting time $dt$ ( $dt \in (0, \Delta T]$ ) |
| **4:** | Reset $\Omega_{tm} = \Omega_{cm}$ |
| **5:** | Inversely propagate $Ele_{cm}, Ele_{tm}$ to $Ele_{c0}, Ele_{t0}$ (the propagation time is $dt$ ) |
| **6:** | Optimize this sample ( $Ele_{c0}, Ele_{t0}$ and $\Delta T$ ) and put into the database pool 2 |
| Algorithm 3 Pseudocode of the process to generate an "RAAN-separating" transfer | |
| **1:** | Randomly produce $Ele_{c0}$, $Ele_{t0}$ and $\Delta T$ within each range determined by the working situation of the problem |
| **2:** | Calculate $\dot{\Omega}_c$ and $\dot{\Omega}_t$ |
| **3:** | Randomly produce an initial RAAN difference $\Delta\Omega_0$ ( $\Delta\Omega_0 \in (0, d_2]$ ) |
| **4:** | **if** $\dot{\Omega}_c < \dot{\Omega}_t$ |
| | Reset $\Omega_{t0} = \Omega_{c0} + \Delta\Omega_0$ |
| | **else** |
| | Reset $\Omega_{t0} = \Omega_{c0} - \Delta\Omega_0$ |
| | **end if** |
| **5:** | Optimize this sample ( $Ele_{c0}, Ele_{t0}$ and $\Delta T$ ) and put into the database pool 3 |

## C. DNN models and network training method

A DNN is made up of large numbers of simple, highly interconnected processing nodes. Each node takes one or more inputs from other nodes and produces an output by applying an activation function over the weighted sum of these inputs. In this study, multilayer perceptron (MLP) is selected as the architecture for the regression DNNs. The activation of a node in MLP is determined by the summation of all the weighted inputs, which can be expressed as

$$x_j = f(\sum_{i=1}^{N} w_{ij} x_i + b_j) \, , \tag{9}$$

where $x_j$ is the output of node $j$ in the current layer, $x_i$ is the output of node $i$ in the previous layer, $w_{ij}$ refers to the weight of the connection from node $i$ to node $j$, $b_j$ denotes the variable bias of node $j$, $N$ is the total number of nodes in the previous layer, and $f$ is the activation function. A Leaky Rectified Linear Unit (ReLU) [28] and a linear function are selected as the hidden-layer and output-layer activation functions, respectively.

Network training can be regarded as a process to adjust the weight vectors epoch by epoch, and the aim is to minimize the loss function. The mean squared error function $F_r$ are used as the loss function for the regression DNNs, and is expressed as

$$F_r = \frac{1}{b} \sum_{i=1}^{b} (o_p(i) - o_m(i))^2 ,\qquad (10)$$

where $b$ is the batch size and $o_p(i)$ is the predicted output of the network. $o_m(i)$ is the optimal velocity increments of a transfer. Cross-validation is applied in each epoch, and 90% of the data are used as training samples while the remaining 10% are used for validation. All the three regression DNNs are trained until convergence with mini-batch gradient decent and a batch size of $b = 32$. The adaptive moment (Adam) [29] technique is used to optimize the parameters of the networks. Keras [30] combined with TensorFlow [31] is applied to train the networks, where TensorFlow is the backend of Keras.

## VI.  Simulations

The missions in GTOC-9 [2] are selected as examples to demonstrate the proposed DNN-based method for approximating optimal transfers.

## A. Generation of the database

Following the configuration in GTOC-9, $\Delta T_{\max}$ is set to 30 days. The orbit elements of the departure and rendezvous debris are all within the ranges shown in Table 4. The acceptable terminal errors are set to 1 m and 0.01 m/s in position and velocity respectively.

**Table 4 Ranges of the orbit elements for both the departure and rendezvous debris**

| $a$ , km | $e$ | $i$ , deg | $\Omega$ , deg | $\omega$ , deg | $f$ , deg |
|---|---|---|---|---|---|
| 6900~7300 | 0~0.02 | 96~101 | 0~360 | 0~360 | 0~360 |

The parameters $d_1$ and $d_2$ in Algorithms 1-3 are both set to 10 deg. Large numbers of transfers are obtained using the database generation method. The samples belonging to each transfer type are generated individually using Algorithms 1-3. One thousand transfers are randomly selected as the test samples for each problem, and the remaining ones are employed as the training samples.

## B. Selection of the learning features

An appropriate selection of the learning features is important because the lack of the relevant features and the interference of the redundant features both reduce the approximating performance [32]. Domain knowledge suggests that the optimal velocity increments of a perturbed long-duration impulsive transfer should relate to its initial and final states. The possible appropriate features for approximating optimal transfers are listed in Table 5, where the first six parameters are required features and others are alternative ones. The number of training samples is set to 5000, and a two-hidden-layer network with 30 nodes is first applied to compare the approximating performance of different feature combinations.

**Table 5 Possible features for approximating optimal transfers**

| Feature | Description |
|---------|-------------|
| $a_c$, $a_t$ | Semi-major axes of the departure body and rendezvous target |
| $e_c$, $e_t$ | Eccentricities of the departure body and rendezvous target |
| $i_c$, $i_t$ | Inclinations of the departure body and rendezvous target |
| $\Delta\Omega_{c0t0}$ | Difference between initial RAAN of the departure body and initial RAAN of the rendezvous target |
| $\Delta\Omega_{cftf}$ | Difference between final RAAN of the departure body and final RAAN of the rendezvous target |
| $\Delta\Omega_{c0tf}$ | Difference between initial RAAN of the departure body and final RAAN of the rendezvous target |
| $\dot{\Omega}_c$, $\dot{\Omega}_t$ | RAAN variation rates of the departure body and rendezvous target |
| $\Delta\varphi_0$, $\Delta\varphi_f$ | Initial and final phase differences between the departure body and rendezvous target |
| $\Delta T$ | Transfer time |

Mean relative error (MRE) is used as the evaluation criterion and is defined as

$$\varepsilon_{\mathrm{MRE}} = \frac{1}{N} \cdot \sum_{i=1}^{N} \frac{|\Delta V_{Esti}^{i} - \Delta V_{Opti}^{i}|}{\Delta V_{Opti}^{i}}, \qquad (11)$$

where $N$ refers to the number of testing samples; $\Delta V_{Esti}^{i}$ and $\Delta V_{Opti}^{i}$ denote the estimated and optimized optimal velocity increments of the $i$th transfer, respectively. Tables 6-8 list the MREs of all the tested groups. The MREs of Groups 5 and 6 in Tables 6-8 show that $\dot{\Omega}_c$, $\dot{\Omega}_t$ and $\Delta\varphi_0$, $\Delta\varphi_f$ contribute little to the improvement of the approximating performance for all the three types of transfers. These results strongly verify that the optimal velocity increments of long-duration impulsive transfers indeed relate little to $\Delta\varphi_0$ and $\Delta\varphi_f$. The comparisons of Groups 1, 2, 3, 4 and 7 in Tables 6 and 8 show that $\Delta\Omega_{c0t0}$, $\Delta\Omega_{cftf}$ and $\Delta T$ can help improve the approximating performance for both "RAAN-closing" and "RAAN-separating" transfers, while $\Delta\Omega_{c0tf}$ is helpful only for approximating "RAAN-separating" ones. The comparison of Groups 1-7 in Table 7 verifies the correctness of the first conclusion in Sec. IV and indicates that the required features in Table 5 are sufficient to well estimate the optimal velocity increments of "RAAN-intersecting" transfers. Consequently, the combinations of Group 7 in Table 6, Group 1 in Table 7 and Group 7 in Table 8 are selected as the learning features, and they are fixed in the following simulations. The results in these three Tables also reflect the difference of the three regression problems and verify the reasonability and necessity of the division for perturbed long-duration impulsive transfers.

**Table 6 MREs of the estimation of "RAAN-closing" transfers using different features**

| Group | Feature combination | MRE |
|---|---|---|
| 1 | $a_c$, $a_t+e_c$, $e_t+i_c$, $i_t$ | 46.22% |
| 2 | $a_c$, $a_t+e_c$, $e_t+i_c$, $i_t+\Delta\Omega_{c0t0}$ | 28.49% |
| 3 | $a_c$, $a_t+e_c$, $e_t+i_c$, $i_t+\Delta\Omega_{c0t0}+\Delta\Omega_{ctf}$ | 12.54% |
| 4 | $a_c$, $a_t+e_c$, $e_t+i_c$, $i_t+\Delta\Omega_{c0t0}+\Delta\Omega_{ctf}+\Delta\Omega_{c0tf}$ | 13.15% |
| 5 | $a_c$, $a_t+e_c$, $e_t+i_c$, $i_t+\Delta\Omega_{c0t0}+\Delta\Omega_{ctf}+\dot\Omega_c$, $\dot\Omega_t$ | 12.51% |
| 6 | $a_c$, $a_t+e_c$, $e_t+i_c$, $i_t+\Delta\Omega_{c0t0}+\Delta\Omega_{ctf}+\Delta\varphi_0$, $\Delta\varphi_f$ | 12.49% |
| 7 | $a_c$, $a_t+e_c$, $e_t+i_c$, $i_t+\Delta\Omega_{c0t0}+\Delta\Omega_{ctf}+\Delta T$ | 5.98% |

**Table 7 MREs of the estimation of "RAAN-intersecting" transfers using different features**

| Group | Feature combination | MRE |
|---|---|---|
| 1 | $a_c$, $a_t+e_c$, $e_t+i_c$, $i_t$ | 5.95% |
| 2 | $a_c$, $a_t+e_c$, $e_t+i_c$, $i_t+\Delta\Omega_{c0t0}$ | 5.94% |
| 3 | $a_c$, $a_t+e_c$, $e_t+i_c$, $i_t+\Delta\Omega_{c0t0}+\Delta\Omega_{ctf}$ | 6.04% |
| 4 | $a_c$, $a_t+e_c$, $e_t+i_c$, $i_t+\Delta\Omega_{c0t0}+\Delta\Omega_{c0tf}$ | 6.09% |
| 5 | $a_c$, $a_t+e_c$, $e_t+i_c$, $i_t+\Delta\Omega_{c0t0}+\dot\Omega_c$, $\dot\Omega_t$ | 5.99% |
| 6 | $a_c$, $a_t+e_c$, $e_t+i_c$, $i_t+\Delta\Omega_{c0t0}+\Delta\varphi_0$, $\Delta\varphi_f$ | 5.97% |
| 7 | $a_c$, $a_t+e_c$, $e_t+i_c$, $i_t+\Delta\Omega_{c0t0}+\Delta T$ | 6.00% |

**Table 8 MREs of the estimation of "RAAN-separating" transfers using different features**

| Group | Feature combination | MRE |
|---|---|---|
| 1 | $a_c$, $a_t+e_c$, $e_t+i_c$, $i_t$ | 44.54% |
| 2 | $a_c$, $a_t+e_c$, $e_t+i_c$, $i_t+\Delta\Omega_{c0t0}$ | 14.64% |
| 3 | $a_c$, $a_t+e_c$, $e_t+i_c$, $i_t+\Delta\Omega_{c0t0}+\Delta\Omega_{ctf}$ | 11.25% |
| 4 | $a_c$, $a_t+e_c$, $e_t+i_c$, $i_t+\Delta\Omega_{c0t0}+\Delta\Omega_{ctf}+\Delta\Omega_{c0tf}$ | 8.54% |
| 5 | $a_c$, $a_t+e_c$, $e_t+i_c$, $i_t+\Delta\Omega_{c0t0}+\Delta\Omega_{ctf}+\Delta\Omega_{c0tf}+\dot\Omega_c$, $\dot\Omega_t$ | 8.53% |
| 6 | $a_c$, $a_t+e_c$, $e_t+i_c$, $i_t+\Delta\Omega_{c0t0}+\Delta\Omega_{ctf}+\Delta\Omega_{c0tf}+\Delta\varphi_0$, $\Delta\varphi_f$ | 8.52% |
| 7 | $a_c$, $a_t+e_c$, $e_t+i_c$, $i_t+\Delta\Omega_{c0t0}+\Delta\Omega_{ctf}+\Delta\Omega_{c0tf}+\Delta T$ | 5.51% |

## C. Determination of the network and training data scales

An appropriate scale of the network is necessary to avoid underfitting and overfitting. Different numbers of hidden layers and nodes are thus tested to determine the most appropriate scales for the three regression problems.

Figures 5-7 illustrate the MREs of the estimations with the number of hidden layers varying from two to four and the number of nodes in each layer varying from 10 to 100. It is shown that a network scale of $2 \times 60$ is appropriate for approximating "RAAN-intersecting" transfers. However, larger scales are required for the other two types of transfers, where $3 \times 60$ and $3 \times 70$ should be the best choices, respectively. Figure 8 illustrates the decreases in the MREs as the number of training samples increases from $2 \times 10^3$ to $10^5$. The results show that the MREs of the estimations can be decreased to less than 3% for all the three types of transfers with sufficient training samples. $5 \times 10^4$ samples is enough to converge to the near-best approximating performance for "RAAN-intersecting" transfers, while at least $7 \times 10^4$ and $8 \times 10^4$ samples are required for "RAAN-closing" and "RAAN-separating" ones, respectively. The requirement of larger network and training data scales indicates that approximating an "RAAN-closing" or "RAAN-separating" transfer is more difficult than approximating an "RAAN-intersecting" one.

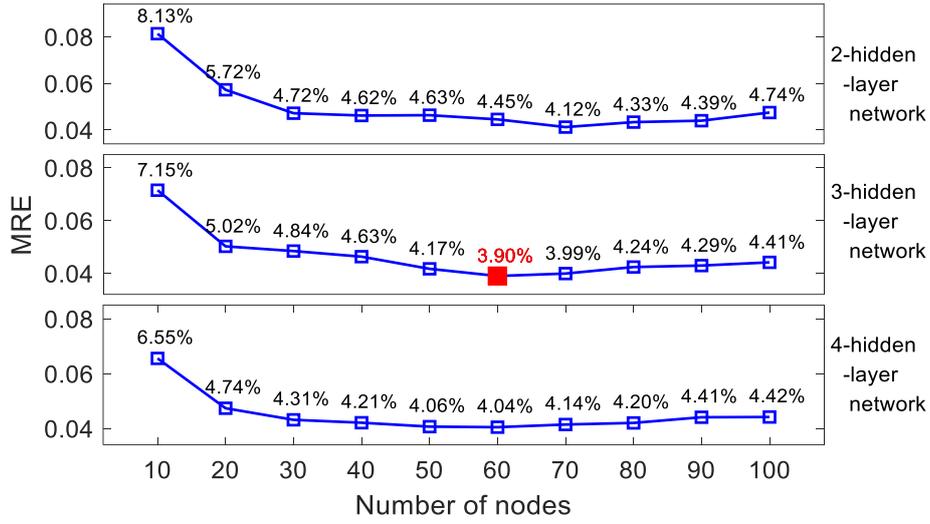

**Fig. 5 MREs of the estimation of "RAAN-closing" transfers for networks with different numbers of hidden layers and nodes**

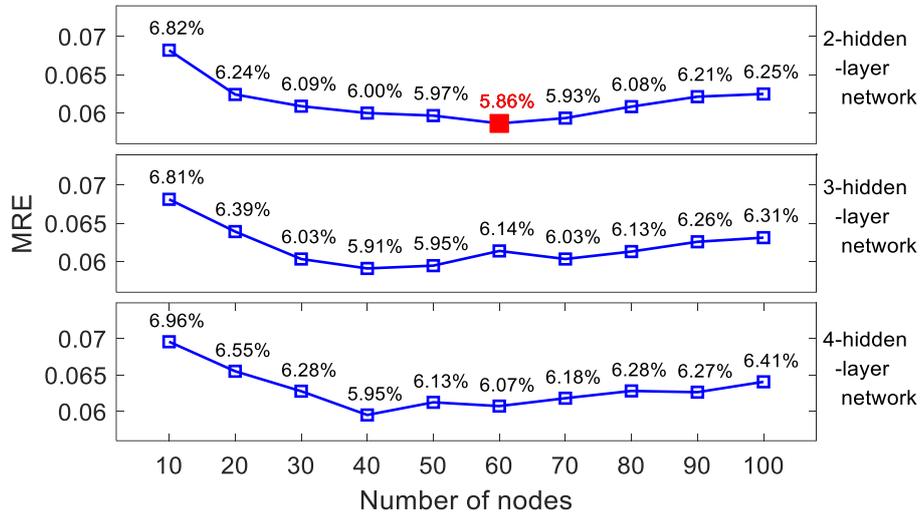

**Fig. 6 MREs of the estimation of "RAAN- intersecting" transfers for networks with different numbers of hidden layers and nodes**

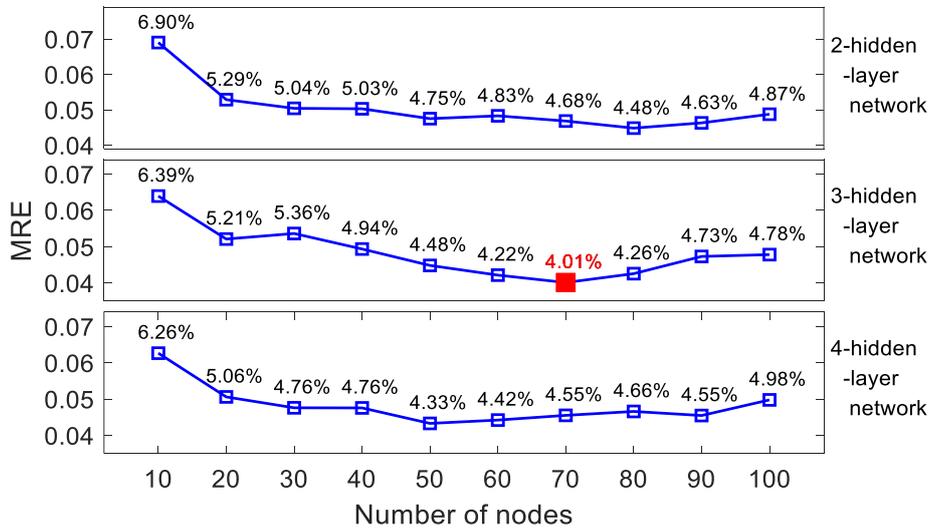

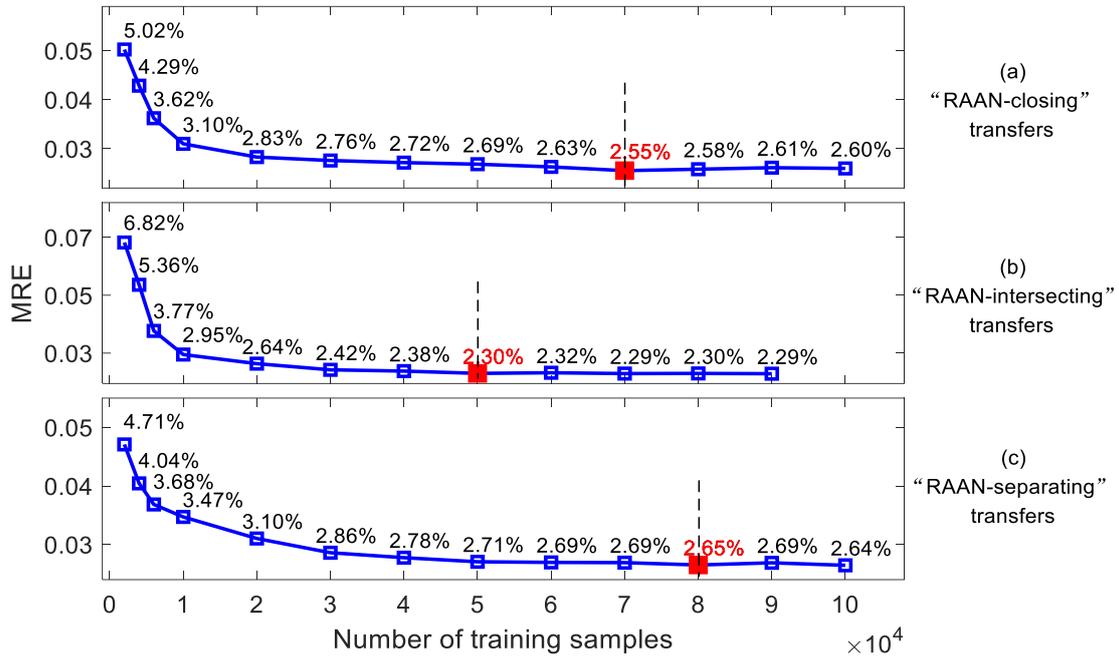

**Fig. 8 MREs of the estimation with different training data scales for the three types of transfers**

Due to the unavoidable errors between the training samples and the corresponding true optimal transfers, and also because of the possible systematic noise of the learning problems, further enlarging the training data scale makes no significant contribution to improvement in the approximating performance. Nevertheless, such an approximating performance is competent for application to the preliminary design of multitarget rendezvous missions.

Moreover, we notice that some participants of GTOC-9 including the champion team Jet Propulsion Laboratory (JPL) have created a database for quick lookups for approximating the optimal transfers between any two pieces of debris [3, 6]. Billions of debris-to-debris transfers are required to approximately cover the whole search space. While only tens of thousands of training samples are enough to approximate the optimal transfers with an acceptable estimation error by means of the DNN-based method. This comparison strongly shows the advantage of the DNN-based method for application to real-world problems.

**D.  Comparison with analytical estimation methods**

Two typical analytical methods are compared with the DNN-based method. The first method proposed by the National University of Defense Technology (NUDT) team is mainly based on the Gauss form of variational equations [33], and the second one proposed by the Xi'an Satellite Control Center (XSCC) team is similar to the Edelbaum's method for approximating high-thrust transfers [34]. Detailed procedures of the two methods can be found in [4] and [5]. For convenience, they are named as NUDT method and XSCC method in this paper.

The estimation results of the test samples obtained by these methods are listed in Table 9. It can be seen from Table 9 that both NUDT method and XSCC method perform much worse than the DNN-based method on all the three types of transfers. Due to the lack of the consideration of necessary features such as the transfer time, the MREs of the two analytical methods can be as high as over 20% when approximating "RAAN-closing" and "RAAN-separating" transfers. Figure 9 further visualizes the statistics of the estimation error for the three methods. We can find that the distributions of the estimation error obtained by the two analytical methods are much wider than those obtained by the DNN-based method, and the centers of the distributions both deviate from 0. These results indicate that there are systematic errors using either NUDT method or XSCC method and the random errors for the two analytical methods are both larger.

**Table 9 MREs of the estimation for the two analytical methods and DNN-based method**

|  | "RAAN-closing" transfers | "RAAN-intersecting" transfers | "RAAN-separating" transfers |
|---|---|---|---|
| NUDT method [4] | 20.37% | 17.93% | 20.82% |
| XSCC method [5] | 21.96% | 16.28% | 22.43% |
| DNN-based method | 2.56% | 2.29% | 2.64% |

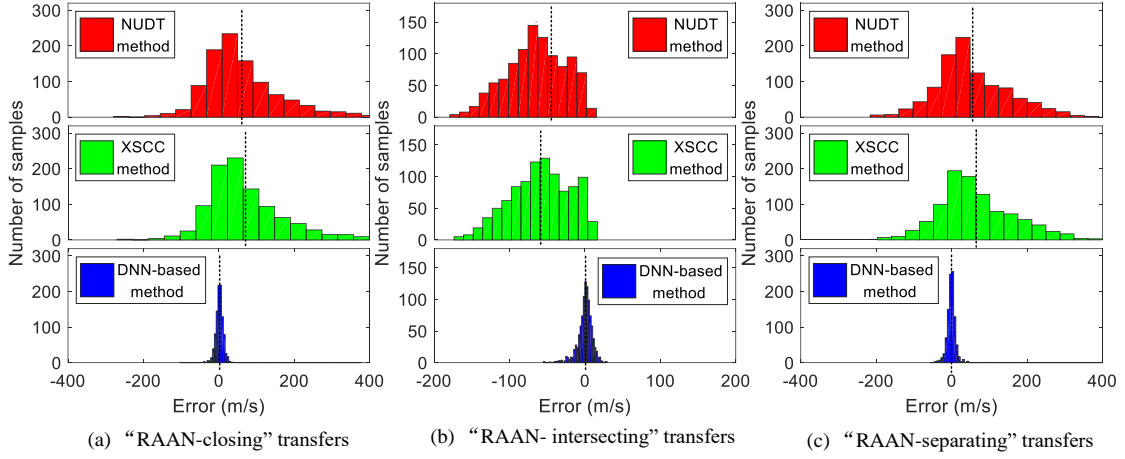

(a) "RAAN-closing" transfers    (b) "RAAN- intersecting" transfers    (c) "RAAN-separating" transfers

**Fig. 9 Error distributions of the tested samples obtained by the two analytical methods and DNN-based method**

## E. Verification on debris transfer chains

To verify the efficiency of the DNN-based method for real-world applications, the fast approximation of successive transfers is further studied. Two transfer chains that are achieved by the JPL team in GTOC-9 are tested [3], where a total of 14 and 12 pieces of debris are contained in Chain 1 and 2, respectively. The debris IDs, as well as the departure and rendezvous time of each transfer are listed in Table 10. The epoch data of the debris can be accessed on GTOC Portal [35].

**Table 10 Transfer chains in GTOC-9 obtained by JPL [3]**

| Transfer Chain 1 | | | Transfer Chain 2 | | |
|---|---|---|---|---|---|
| Debris ID | Rendezvous time, MJD 2000 | Departure time, MJD 2000 | Debris ID | Rendezvous time, MJD 2000 | Departure time, MJD 2000 |
| 23 | — | 23557.18 | 33 | — | 25951.06 |
| 55 | 23587.04 | 23592.04 | 68 | 25973.75 | 25979.26 |
| 79 | 23617.02 | 23622.06 | 116 | 25983.50 | 25989.03 |
| 113 | 23644.48 | 23649.49 | 106 | 26013.50 | 26019.03 |
| 25 | 23674.48 | 23679.49 | 14 | 26043.49 | 26049.02 |

| 20 | 23679.78 | 23684.81 | 52 | 26073.49 | 26079.04 |
|---|---|---|---|---|---|
| 27 | 23695.44 | 23700.44 | 120 | 26103.48 | 26109.02 |
| 117 | 23725.44 | 23730.44 | 80 | 26133.48 | 26139.01 |
| 121 | 23733.14 | 23738.14 | 16 | 26163.47 | 26169.01 |
| 50 | 23739.65 | 23744.68 | 94 | 26193.47 | 26199.02 |
| 95 | 23746.09 | 23751.12 | 83 | 26217.56 | 26223.08 |
| 102 | 23775.79 | 23780.83 | 89 | 26232.30 | — |
| 38 | 23805.14 | 23810.18 | | | |
| 97 | 23816.04 | — | | | |

The optimized and estimated results of the two chains are illustrated in Figures 10 and 11. Note that the estimation accuracy of the whole chain but not that of each transfer is more focused when approximating successive transfers. It can be seen from Figures 10 and 11 that the total velocity increments estimated by the DNN-based method are both very close to the results optimized by the JPL team (the final errors are approximately 9 m/s and 8 m/s for Chain 1 and Chain 2, respectively). These results indicate that the DNN-based method can be well applied to real-world problems with excellent approximating performance.

More detailedly, we find that the estimated and optimized velocity increments of each transfer in Chain 1 are not fitted so well as those in Chain 2, especially for the first sixth transfers (the accumulative error has come to approximately 70 m/s). Nevertheless, from Figure 9 we know that there is no systematic error but only a random error using the DNN-based method. Even if the accumulative error of the first sixth transfers is a little large, it can be offset by the following ones and finally reaches a small value of less than 10 m/s. The longer the chain is, the less the final estimated result will be influenced by the random error. However, the accumulative error will keep increasing and finally reach an unacceptable value if applying analytical methods because of the existence of the systematic error. This comparison further shows the significant advantage of the DNN-based method for approximating successive transfers.

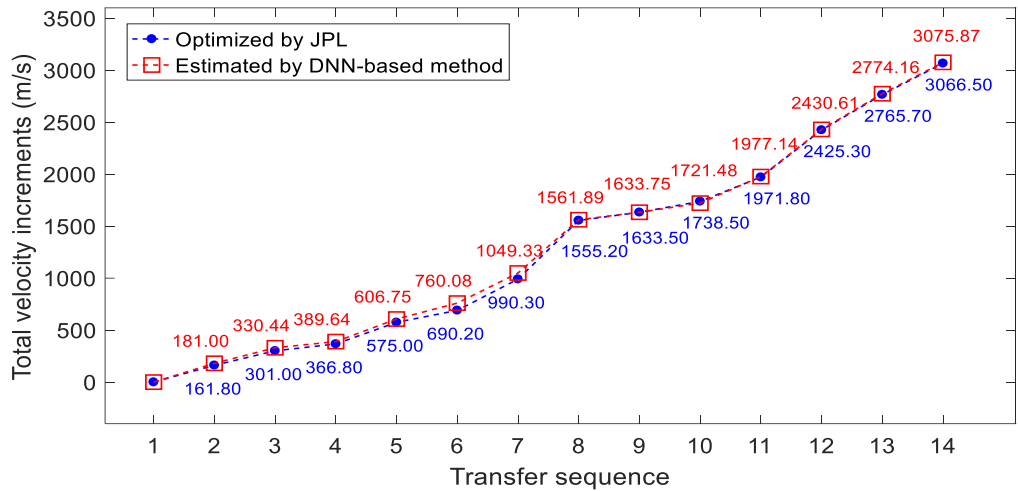

**Fig. 10 Optimized and estimated total velocity increments of Chain 1**

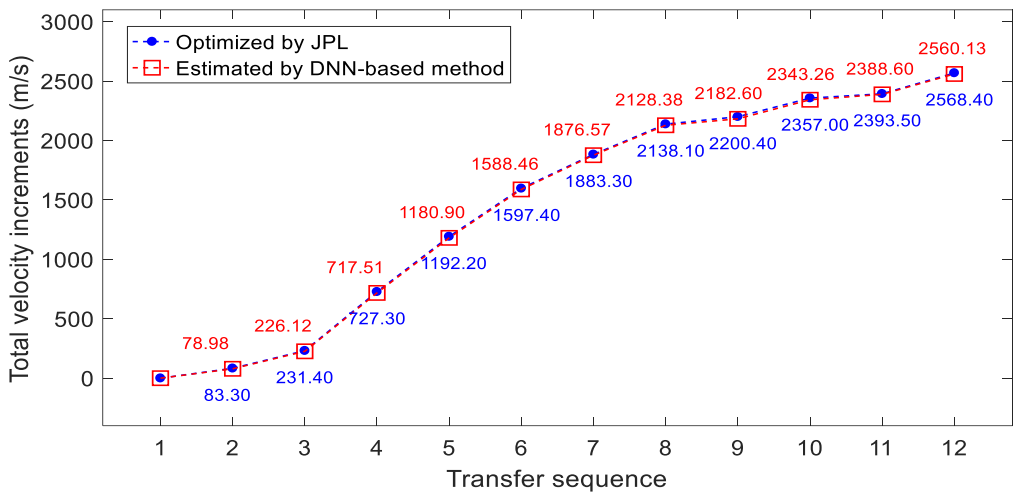

**Fig. 11 Optimized and estimated total velocity increments of Chain 2**

## VII.    Conclusions

Fast approximation of optimal perturbed long-duration impulsive transfers is studied. The relationship between the optimal velocity increments of this kind of transfer and the initial RAAN difference between the departure body and the rendezvous target, as well as the transfer time is investigated, and the results show that this kind of transfer should be divided into three types, which are "RAAN-closing" transfers, "RAAN-intersecting" transfers and "RAAN-separating" transfers, respectively. A DNN-based method for quickly approximating optimal transfers is proposed, where three regression DNNs are designed to estimate the optimal velocity increments. Each DNN corresponds to a type of transfers, and the DNNs are trained individually based on the corresponding database. The most appropriate learning features and network scales for approximating each type of transfers are determined. The comparison with analytical estimation methods shows the superiority and reliability of the DNN-based method for quickly approximating this kind of transfer. The study on two debris chains with more than ten transfers further reveals the advantage of the DNN-based method for estimating the optimal velocity increments for successive transfers.

## Acknowledgments

This work was supported by the Natural Science Foundation of China (11572345).